# Marvelous cancellations: T.J. Stieltjes' letters concerning the zeta function

Juan Miguel Marín

*(jmarin@mail.harvard.edu)*

**Abstract** This article introduces and translates letters from T.J Stieltjes (1856-1894) to C. Hermite (1822-1891) regarding Stieltjes' published notes claiming to have solved Bernhard Riemann's conjecture "it is very probable that all the roots [of the zeta function] are real." These letters illuminate what was at stake in the transition between the "non-rigorous" method associated with the 18th century and the "rigorous" approach attributed to the late 19th. They reveal not only obsolete assumptions that no longer remain with us but also those unexamined ones that still do. Becoming aware of the issues at stake can help us understand the historical, philosophical and scientific background from which today's number theory emerged.

*MSC:* 01A55

*Keywords*:  Correspondence; Stieltjes; Hermite

## 1. Introduction

Dutch mathematician T.J. Stieltjes' claim  to have solved Riemann's hypothesis has received attention in the scholarly literature , especially after the 1985 proof of the falsity of the Merten's conjecture, a stronger version of Stieltjes' 1885 claim to have proven, in today's notation, that $M(n)/\sqrt{n}$ is bounded. [Borwein, 2009, 69] Number theorists today are skeptical of Stieltjes' claim that he had a proof of an equivalent of Riemann's hypothesis, and with good reason. And yet, revisiting his attempt at proof, something that has never been done, can yield insights on the history of understanding the problem.  His letters offer a little studied resource on this first attempt to solve Riemann's conjecture, an attempt that, while lacking rigor by contemporary standards, still has not been entirely rejected. Studying Stieltjes' letters can still shed light, not only on the assumptions of his 19th century mathematical contemporaries but also on those assumptions we have inherited from them.

### 1.1    Intro [1]

On July 13, 1885; *Comptes Rendus* published Stieltjes' note *Sur une function uniforme*. It opens as follows:

The analytical character of the function $\zeta(z)$ which is defined for the values of $z$, with real part greater than unity, by the series

$$1 + \frac{1}{2^z} + \frac{1}{3^z} + \frac{1}{4^z} + \cdots$$

has been revealed entirely by Riemann who showed that

$$\zeta(z) - \frac{1}{z-1}$$

is holomorphic throughout the entire plane.

The zeroes of the function $\zeta(z)$ are on one hand

$$-2, -4, -6, -8 \ldots$$

and on the other, there are an infinity of zeros , all of them imaginary, the remaining real part within 0 and 1.

Riemann announced as very probable that all these imaginary roots are of the form $\frac{1}{2} + \alpha i$, $\alpha$ being real. I managed to prove rigorously Riemann's claim and will describe here the way that led me to this result.

Following Euler, we have that

$$1 : \zeta(z) = \prod \left(1 - \frac{1}{p^z}\right),$$

$p$ representing all prime numbers, that is

$$1 : \zeta(z) = 1 - \frac{1}{2^z} - \frac{1}{3^z} - \frac{1}{5^z} + \frac{1}{6^z} - \frac{1}{7^z} + \frac{1}{10^z} - \cdots$$

It was an in-depth study of the second series above which led to the desired goal.

[Stieltjes, 1885, 153]

As Narciewicz shows in his excellent discussion, Stieltjes claimed that the above series converges when the real part of $s$ is greater than $\frac{1}{2}$, , a claim that implies the truth of Riemann's conjecture. Today we know that the convergence of

$$\sum_{n=1}^{\infty} \frac{\mu(n)}{n^s}; Re(s) > \frac{1}{2}$$

is not only necessary but sufficient for the Riemann hypothesis.(Borwein, 2008)(Titchmarsh, 1987) . Now, we notice that in Stieltjes' notation the variable he uses is, not Riemann's $s$ but $z$. The variable betrays its origin.

Stieltjes obtained the above from Euler's *Analysis of the Infinities*. Euler writes the above as

$$1: \zeta(z) = 1 - \frac{1}{2^z} - \frac{1}{3^z} - \frac{1}{5^z} + \frac{1}{6^z} - \frac{1}{7^z} + \frac{1}{10^z} - \cdots$$

[Euler, 1999]

which then Stietjes' develops, as we can discover by studying the letters. These posthumously letters can provide us with much insight into Stieltjes' published work.

When Hermite asked him to elaborate on how he proved the above claim Stieltjes answered in his July 11, 1885 letter:

> When considering the expression obtained by Riemann for $\xi(t)$ one finds that it has [some] real roots; but I have fruitlessly tried to deduce this expression, by definite integral. That it has *all* its roots real, and I became desperate in trying to prove such a doubtful proposition, when I perceived that one obtains this proposition by modifying slightly Riemann's reasoning for obtaining the zeroes of $\zeta(s)$ within this mysterious strip where the real part of $s$ lies between 0 and 1. In fact, if, instead of $1: \zeta(s) = \prod \left(1 - \frac{1}{p^s}\right)$, one considers $1: \zeta(s) = 1 - \frac{1}{2^s} - \frac{1}{3^s} - \frac{1}{5^s} + \frac{1}{6^s} + \cdots = \sum_1^{\infty} \frac{f(n)}{n^s}$, there is a major difference between the infinite product and the series; the latter is convergent for $s > \frac{1}{2}$; while for the former we ought to assume that $s > 1$.

"Here's how I proved it: The function $f(n)$ is equal to zero when $n$ is divisible by a square and for the other values of $n$, equal to $(-1)^k$, $k$ being the number of prime factors of $n$. Now, I found that for the sum

$$g(n) = f(1) + f(2) + \cdots + f(n),$$

the terms $\pm 1$ cancel each other (*se compensent assez bien*) for $\frac{g(n)}{\sqrt{n}}$ always remains between two fixed limits, regardless how large $n$ may be (probably one could take $+1$ and $-1$). From this it follows that $s > \frac{1}{2}$,

$$\lim \frac{g(n)}{n^s} = 0 \qquad (n = \infty)$$

and similarly, with $|g(n)|$ designating the absolute value of $g(n)$, the series $\sum \frac{|g(n)|}{n^{1+s}}$ converges.

Perhaps Stieltjes here follows his century's "bible" of mathematics, Cauchy's *Cours d' Analyse*.

In the Cours' Note II, Theorem XV Cauchy works with the quadratic or "square root" mean as follows:

*Let $a, a', a'', \ldots$ be any $n$ real quantities. If these quantities are not equal to each other,*

*then the numerical value of the sum*

$$a + a' + a'' + \cdots$$

*is less than the product*

$$\sqrt{n}\sqrt{a^2 + a'^2 + a''^2 + \cdots}$$

*so that we have*

$$val.num. \; (a + a' + a'' + \cdots) < \sqrt{n}\sqrt{a^2 + a'^2 + a''^2 + \cdots}$$

[Cauchy, 2009, 302]

The last part is equivalent to $|a + a' + a'' + \cdots| < \sqrt{n}\sqrt{a^2 + a'^2 + a''^2 + \cdots}$ . So if we let $a = (-1)^k$, $k$ being the number of distinctive prime factors, we have

$$|a_1 + a_2 + \cdots + a_n| \leq \sqrt{n}\sqrt{a_1^2 + a_2^2 + \cdots + a_n^2}$$

$$|a_1 + a_2 + \cdots + a_n| \leq \sqrt{n}\sqrt{(-1_1)^{2k} + (-1_2)^{2k} + \cdots + (-1_n)^{2k}}$$

$$|a_1 + a_2 + \cdots + a_n| < \sqrt{n}\sqrt{(-1_1)^2 + (-1_2)^2 + \cdots + (-1_n)^2}$$

So Stieltjes may have been thinking of $a_n = \mu_n$

$$\sum_{i=1}^{n} a_n \leq \sqrt{n}\sqrt{\sum_{i=1}^{n} 1}$$

$$\sum_{i=1}^{n} a_n \leq \sqrt{n}\sqrt{\zeta(0)}$$

$$\sum_{i=1}^{n} a_n \leq \sqrt{n}\sqrt{-\frac{1}{2}}$$

The latter can be explained via zeta function regularization, which today primarily appears in the work of mathematical physicists. Stieltjes himself gave an "electrostatic interpretation" of the zeros of Jacobi, Laguerre, and Hermite polynomials (Stieltjes, 13). We can therefore wonder whether a related "physical" interpretation would have been in his mind, one as controversial as Riemann's use of the "physical" Dirichlet's principle.

**Conclusion**

Stieltjes' letters illuminate what was at stake in the transition between the "non-rigorous" method associated with the 18[th] century and the "rigorous" approach attributed to the late 19[th]. They reveal not only obsolete assumptions that no longer remain with us but also those unexamined ones that still do. Becoming aware of the issues at stake can help us understand the historical, philosophical and scientific background from which today's number theory emerged.

In answer to Stieltjes' letters, Hermite answered that his mathematics were indeed beautiful, only one thing remained to be done, one noted by Mittag-Leffler.

Mittag-Leffler could not go beyond a point; he cannot see how you established that the series

$$1 - \frac{1}{2^s} - \frac{1}{3^s} - \frac{1}{5^s} + \frac{1}{6^s} - \cdots$$

converges when the real part of $s$ is greater than $\frac{1}{2}$. On his behalf, I beg you to write him and free him from this embarrassment. You would him much good if you send him the clarifications, which he awaits, as he says, with impatience …

We all still are.

## 3. Appendix: Stieltjes' Notes and Letters to Hermite (1885-1887)

### Notes in Compte Rendus

*On a uniform function.* Note by T.J. Stieltjes presented by C. Hermite.

The analytical character of the function $\zeta(z)$ which is defined for the values of $z$, with real part greater than unity, by the series

$$1 + \frac{1}{2^z} + \frac{1}{3^z} + \frac{1}{4^z} + \cdots$$

has been entirely revealed by Riemann who has shown that

$$\zeta(z) - \frac{1}{z-1}$$

is holomorphic throughout the entire plane.

The zeroes of the function $\zeta(z)$ are on one hand

$$-2, -4, -6, -8 \ldots$$

and on the other, there are an infinity of zeros , all of them imaginary, the remaining real part within 0 and 1.

Riemann announced as very probable that all these imaginary roots are of the form $\frac{1}{2} + \alpha i, \alpha$ being real. I managed to prove rigorously Riemann's claim and will describe here the way that led me to this result.

Following Euler, we have that

$$1 : \zeta(z) = \prod \left( 1 - \frac{1}{p^z} \right),$$

$p$ representing all prime numbers, that is

$$1 : \zeta(z) = 1 - \frac{1}{2^z} - \frac{1}{3^z} - \frac{1}{5^z} + \frac{1}{6^z} - \frac{1}{7^z} + \frac{1}{10^z} - \cdots$$

It was an in-depth study of the second series above which led to the desired goal. In fact, one can prove that this series converges and defines an analytical function as long as the real part of $s$ is greater than $\frac{1}{2}$ . It is clear from the above that $\zeta(z)$ does not vanish for any value of $z$ where the real part of $z$ is greater than $\frac{1}{2}$. But the equation $\zeta(z) = 0$ cannot admit any more imaginary roots where the real part is less than $\frac{1}{2}$. In fact, if there were such a root $z = z_1$ one would also have $\zeta(1 - z_1) = 0$, as shown by the relationship between $\zeta(z)$ and $\zeta(1 - z)$ established by Riemann. Now, the real part of $(1 - z_1)$ is greater than $\frac{1}{2}$. Therefore, *all imaginary roots of $\zeta(z) = 0$ are of the form $\frac{1}{2} + \alpha i$ with $\alpha$ real.*

*Concerning an asymptotic law in number theory.* Note by T.J. Stieltjes presented by C. Hermite.

The theorem announced in the *Comptes Rendus*, p. 53, that the series $(A)$

$$1 - \frac{1}{2^s} - \frac{1}{3^s} - \frac{1}{5^s} + \frac{1}{6^s} - \cdots$$

obtained by the development of the infinite product

$$\prod \left(1 - \frac{1}{p^s}\right),$$

converges for $s > \frac{1}{2}$ leads to an important consequence related to Tchebychef's function $\theta(x) = $ the sum of the logarithms of prime numbers less than $x$.

Let $f(n)$ be the number of divisors of $n$.we recall the result due to Dirichlet that

$$\frac{f(1) + f(2) + \cdots + f(m) - n \log n - (2C - 1)n}{\sqrt{n}}$$

remains within two fixed limits, with $C$ being the Eulerian constant. One can then easily conclude that the series $(B)$

$$\sum_1^\infty \frac{f(n) - \log n - 2C}{n^s}$$

converges for $s > \frac{1}{2}$

Here are now two theorems related to the series of the form $\sum_1^\infty \frac{\lambda(n)}{n^s}$. Both theorems are necessary for our discussion.

**THEOREM I.** *Whenever the series* $\sum_1^\infty \frac{\lambda(n)}{n^s}, s > 0$, is convergent, we have

$$\lim \frac{\lambda(1) + \lambda(2) + \cdots + \lambda(n)}{n^s} = 0 \qquad (n = \infty).$$

**THEOREM II.** *When the two series*

$$\sum_1^\infty \frac{\lambda(n)}{n^s}, \sum_1^\infty \frac{\mu(n)}{n^s}$$

*are convergent for* $s = \alpha > 0$ *and the series*

$$\sum_1^\infty \frac{|\lambda(n)|}{n^s}, \sum_1^\infty \frac{|\mu(n)|}{n^s}$$

are convergent for $s = \alpha + \frac{1}{2}\beta$, then the series obtained by multiplying the first two

$$\sum_1^\infty \frac{\nu(n)}{n^s},$$

where

$$\nu(n) = \sum \lambda(d)\mu\left(\frac{n}{d}\right),$$

$d$ representing all the divisors of $n$, is convergent for $s = \alpha + \frac{1}{2}\beta$.

Substituting in series $(A)$ and $(B)$ each term for its absolute value, the new series converge for $s > 1$. From Theorem II we then find that multiplying series $(A)$ and $(B)$ the obtained series converges for $> 3/4$.

      Now, one also obtains

$$\sum_1^\infty \frac{1 - g(n)}{n^s}$$

where

$$g(1) = 2C,$$

and, when $p$ is prime , $g(p^k) = \log p$, while $g(n) = 0$ when $n$ is not of the form $p^k$. From Theorem I we conclude that

$$\lim \frac{n - g(1) - g(2) - \cdots - g(n)}{n^s} = 0 \quad \left(n = \infty, \quad s > \frac{3}{4}\right);$$

but we se clearly that

$$g(1) + g(2) + \cdots + g(n) = 2C + \theta(n) + \theta\left(n^{\frac{1}{2}}\right) + \theta\left(n^{\frac{1}{3}}\right) + \cdots$$

with the result that, letting

$$\theta(n) + \theta\left(n^{\frac{1}{2}}\right) + \theta\left(n^{\frac{1}{3}}\right) + \cdots = n + A_n n^s,$$

we find

$$\lim A_n = 0, (n = \infty).$$

It is then easy to deduce that we also have

$$\theta(n) = n + B_n n^s$$

where

$$\lim B_n = 0$$

for $\left(s > \frac{3}{4}\right)$.

This result leads us to this consequence: regardless of how small a positive number $h$ may be, the number of prime numbers within

$$n \text{ and} (1 + h)n$$

always ends up increasing beyond every limit, whenever $n$ grows indefinitely.

## 4. Letters from Stieltjes to Hermite (1885-1887)

*4.1 Stieltjes to Hermite, July [10?], 1885*

Sir,

By digging deeper into the nature of the function $\zeta(s)$ I have been able to finally see the route to reach the results announced by Riemann. I had however to take a different path than that which he indicated. I find the way by which he said to have obtained the approximate number of roots ofauxiliary $\xi(t) = 0$ to be absolutely incomprehensible. I see no means by which to evaluate this integral. I'm therefore happy to have been able to evade this obstacle and prove the claim announced by Riemann as very probable, that is, all root of $\xi(t) = 0$ are real. Through this, the question returns to the discussion of a real function for real values of the variable and this is feasible, at least, and I will move in this direction to attain our goal.

Still, all this research demands much time; among other things, I have to verify my calculations of the constants $C_1, C_2, \ldots, C_5$ and I plan to join the values of the coefficients $D, D_1, D_2, \ldots, D_5$:

$$\zeta(s) = \frac{1}{s-1} + D - D_1 s + \frac{D_2}{1\cdot 2} s^2 - \frac{D_3}{1\cdot 2\cdot 3} s^3 \ldots, \qquad D = \frac{1}{2}$$

$$D_k = (\log 2)^k + (\log 3)^k + \cdots + [\log(n-1)]^k + \frac{1}{2}(\log n)^k - \int_1^n (\log n)^k \, dn$$

$$(k = 1,2,3,\ldots) \;\; D = \frac{1}{2}\log(2\pi) - 1$$

Since I cannot actively continue this work right now because of other duties, I propose to get some fresh air and leave all this for a few months. Nevertheless, I hope it will not be inconvenient to publish in *Comptes rendus* the included Note, which, it seems, will be of interest to geometers who have already studied Riemann's article. The function $\zeta(s)$ is intimately tied to much arithmetical research concerning asymptotic laws related to prime numbers, etc. For example, although I haven't been able to prove it in a rigorous manner, I have no doubt of the claim that $\Phi(x) - \log\log x$ converges for $x = \infty$ through a finite limit [where the expression is a bit complicated, $\Phi(x)$ indicating the sum $\frac{1}{2} + \frac{1}{3} + \frac{1}{5} + \frac{1}{7} + \cdots$ relative to all prime numbers less than $x$.] Besides, Halphen in *Comptes Rendus* of March 5, 1883, has indicated the involvement in these issues of the $\zeta(s)$ function.

I'm glad that you asked me whether you should do me the favor of inserting in the *Comptes rendus*' the attached Note  I can now correct the proofs myself, though naturally, if you run through my Note and make any necessary corrections, I would accept them with appreciation.

Sir, I respectfuly remain, devotedly yours.

*4.2 Hermite to Stieltjes, July 9, 1885*

Paris July 9, 1885

Sir,

Your beautiful discovery concerning Riemann's claim for equation $\xi(t) = 0$ is of utmost interest for me, both for the great importance of the result of having put off my doubting this claim and for the methods you have employed. Nothing would bring me more pleasure than learning the way by which you have achieved the analytic extension of the product $\prod \left( 1 - \frac{1}{p^s} \right), s > \frac{1}{2}$; this way is outside my reach and I haven't been able to come up with any idea. Next Monday your Note will be presented at the Academy's meeting. I have nothing to change in your extremely clear and correct text except to write $\xi(s)$ instead of $\zeta(s)$ so as to use Riemann's own notation in this work. If you wish, I place at your disposition the possibility of intercalating if convenient the proof of the formula $\xi(s) = \frac{1}{s-1} + G(s)$; when going to correct the proofs Wednesday at the Gauthier Villars printing office, you could add it to your text. [2]

With my most sincere congratulations, I renew, Sir, the assurance of my high esteem for you and my feelings of devotion. [3]

*4.3 Stieltjes to Hermite, July 11, 1885*

*Paris; July 11, 1885*

Sir,

Please receive my sincere thanks for the final revision of your demonstration that $\zeta(s) = \frac{1}{s+1} + G(s)$; this monstrance of your kindness is dear to me. Now, le me add that in my Note I always follow Riemann's notation.

Riemann let $\zeta(s) = \sum \frac{1}{n^s}$. After having found that

$$\prod \left(\frac{s}{2} - 1\right) \pi^{-\frac{s}{2}} \zeta(s)$$

does not change when replacing $s$ for $1 - s$, he considered the function obtained through multiplying by $\frac{1}{2} s(s - 1)$,

$$\prod \left(\frac{s}{2}\right)(s - 1) \pi^{-\frac{s}{2}} \zeta(s)$$

which would have the same property. This amounts to say that letting

$$s = \frac{1}{2} + ti \,, \prod \left(\frac{s}{2}\right)(s - 1)\pi^{-\frac{s}{2}} \zeta(s)$$

will be an even function of $t$ designated by $\xi(t)$.

The expression that he found directly for $\xi(t)$, which makes us see that in effect this function is even, gives us therefore a second proof of the relationship obtained above between $\zeta(s)$ and $\zeta(1 - s)$. In regards the function $\xi(t)$, you can see that it has lost the pole $s = 1$ and the zeroes $s = -2, -4, -6, \dots$ Now, the relationship $1 : \zeta(s) = \prod \left(1 - \frac{1}{p^s}\right)$ shows that $\zeta(s)$ does not have zeroes with real part $s > 1$. The relationship between $\zeta(s)$ and $\zeta(1 - s)$ shows then that, within the part of the plane where the real part of $s$ is negative, $s = -2, -4, -6, \dots$ are the *only zeroes*. The function $\xi(t)$ cannot have zeroes except within the strip where the real part is within $0$ and $1$, or, what amounts to the same, if one has

$$\xi(a + bi) = 0$$

$b$ ought to be found within $-\frac{1}{2}$ and $+\frac{1}{2}$. Now, Riemann said that it is very probable that all the zeroes of function $\xi(t)$ are real $b = 0$. Having set $s = \frac{1}{2} + ti$ amounts to say that all imaginary roots of $\zeta(s)$ are of the form $\frac{1}{2} + ai, a$ real. These are the form, slightly different, that I have expressed Riemann's proposition, not wanting to introduce the function $\xi(t)$, which is not the main object of research, and which is introduced only as auxiliary for the function $\zeta(s)$. At least that is how I see it. It is true that this function $\xi(t)$ gathers in itself all the difficulties if one attemptsto obtain the decomposition in primary factors of

$$(s-1)\zeta(s) = \pi^{-\frac{s}{2}}\frac{1}{\prod\left(\frac{s}{2}\right)}\xi(t)$$

When considering the expression obtained by Riemann for $\xi(t)$ one finds that it has [some] real roots; but I have fruitlessly tried to deduce this expression, by definite integral. That it has *all* its roots real, and I became desperate in trying to prove such a doubtful proposition, when I perceived that one obtains this proposition by modifying slightly Riemann's reasoning for obtaining the zeroes of $\zeta(s)$ within this mysterious strip where the real part of $s$ lies between 0 and 1. In fact, if, instead of1: $\zeta(s) = \prod\left(1-\frac{1}{p^s}\right)$, one considers 1: $\zeta(s) = 1 - \frac{1}{2^s} - \frac{1}{3^s} - \frac{1}{5^s} + \frac{1}{6^s} + \cdots = \sum_1^\infty \frac{f(n)}{n^s}$, there is a major difference between the infinite product and the series; the latter is convergent for $s > \frac{1}{2}$;while for the former we ought to assume that $s > 1$.Here's how I proved it: The function $f(n)$ is equal to zero when $n$ is divisible by a square and for the other values of $n$, equal to $(-1)^k$, $k$ being the number of prime factors of $n$. Now, I found that for the sum

$$g(n) = f(1) + f(2) + \cdots + f(n),$$

the terms $\pm 1$ cancel each other (*se compensent assez bien*) for $\frac{g(n)}{\sqrt{n}}$ always remains between two fixed limits, regardless how large $n$ may be (probably one could take $+1$ and $-1$). From this it follows that $s > \frac{1}{2}$,

$$\lim\frac{g(n)}{n^s} = 0 \qquad (n = \infty)$$

and similarly, with $|g(n)|$ designating the absolute value of $g(n)$, the series $\sum\frac{|g(n)|}{n^{1+s}}$ converges.That which we need to prove is that one can make $\sum_1^{n+m}\frac{f(n)}{n^s}$ as small as one wants to by choosing a convenient $n$. But, with $f(n) = g(n) - g(n-1)$ as auxiliary, this expression becomes equal to

$$\frac{g(n+m)}{(n+m)^s} - \frac{g(n-1)}{(n)^s} + g(n)\left[\frac{1}{n^s} - \frac{1}{(n+1)^s}\right]$$

$$+ g(n+1)\left[\frac{1}{(n+1)^s} - \frac{1}{(n+2)^s}\right] + \cdots$$

$$+ g(n+m-1)\left[\frac{1}{(n+m-1)^s} - \frac{1}{(n+m)^s}\right]$$

But we have

$$\frac{1}{n^s} - \frac{1}{(n+1)^s} = \frac{1}{(n+\theta)^{s+1}} \ (0 < \theta < 1).$$

Therefore

$$\sum_1^{n+m} \frac{f(n)}{n^s} = \frac{g(n+m)}{(n+m)^s} - \frac{g(n-1)}{(n)^s}$$

$$+ \frac{sg(n)}{(n+\theta)^{s+1}} + \frac{sg(n+1)}{(n+1+\theta')^{s+1}} + \cdots + \frac{sg(n+1)}{(n+m-1+\theta)^{s+1}}\Bigg\} = \boldsymbol{R}.$$

Now, series $\sum \frac{|g(n)|}{n^{s+1}}$ being convergent, one can render

$$\frac{|g(n)|}{(n)^{s+1}} + \frac{|g(n+1)|}{(n+1)^{s+1}} + \cdots + \frac{|g(n+m-1)|}{(n+m-1)^{s+1}}$$

as small as one wishes; something similar takes place for $\boldsymbol{R}$ set as

$$\sum_1^{n+m} \frac{f(n)}{n^s} = \frac{g(n+m)}{(n+m)^s} - \frac{g(n-1)}{(n)^s} + \boldsymbol{R}$$

Moreover, terms $\frac{g(n+m)}{(n+m)^s}$ and $\frac{g(n-1)}{(n)^s}$ converge to zero and can be made arbitrarily small .

Therefore the series $\sum \frac{f(n)}{n^s}$ converges for $s > \frac{1}{2}$. I believe that it converges again for real value $s = \frac{1}{2}$, but I haven't been able to prove it. What is certain is that it cannot converge when $s < \frac{1}{2}$ and, $s$ being $< \frac{1}{2}$, it is impossible that $\frac{f(1)+f(2)+\cdots f(n)}{n^s}$ lies between two fixed limits [because one would conclude, as in the above, convergence of $\sum \frac{f(n)}{n^s}$ for values less than $\frac{1}{2}$, which is impossible.] This shows clearly the nature of that proposition on which I base myself, that is,

$$\frac{f(1) + f(2) + \cdots f(n)}{\sqrt{n}}$$

remains within two fixed limits.

You can see that everything depends in an arithmetical research on the sum $f(1) + f(2) + \cdots + f(n)$. My proof is tiresome/arduous; when I resume these researches I will try to simplify it.

Nevertheless, one may already make an idea for oneself about the speed with which $g(n)$ increases (or rather that with which the amplitude of its oscillation increases) for the relation $k = E(\sqrt{n})$,

$$g(n) - g\left(\frac{n}{2}\right) + g\left(\frac{n}{3}\right) + \cdots \pm g\left(\frac{n}{k}\right)$$

$$= -1 + h(k)g(k) - h(n)f(1) - h\left(\frac{n}{2}\right)f(2) - \cdots - h\left(\frac{n}{k}\right)f(k),$$

the function $h(x)$ being equal to 1 or 0 according to whether $E(x)$ is even or odd. Since $g(k)$ is naturally $< k$ in absolute value you can see that

$$g(n) - g\left(\frac{n}{2}\right) + \cdots \pm g\left(\frac{n}{k}\right)$$

is less than $2k + 1$ in absolute value, and the same as $k + 1$ when $k$ is even.

You can see clearly now how this study of $\zeta(s)$ has led me to arithmetical speculations. But forgive my having spoken in my previous letter of the proposition

$$\sum_{q<n} \frac{1}{q} - \log\log n = A \ (n = \infty)$$

which has alredy been proven by Mertens (*Crelle*, t. 78). Mertens has also given the determination of $A$, after having already been considered by M. Tchebychef (*Journ. de Liouville*, 1st series, t. XVIII). Legendre already had obtained it by induction.

Sir, I hope that this letter has not been too long and , above all, I hope to have convinced you that I am not too far from Riemann's notation: it seems it's just a slight misunderstanding.

With deep respect I remain devoutly yours.

*4.4 Hermite to Stieltjes, March 29, 1887.*

*Paris, March 29 1887.*

Dear Sir, Mittag-Leffler has taken as subjects of his lectures at the University of Stockholm Riemann's article on the prime numbers, the object of the excellent study you published in *ComptesRendus* in 1885. The wonderful theorem that the number of roots of the equation $\xi(t) = 0$, with real part within the limits where $T$ is $\frac{T}{2\pi} \log \frac{T}{2\pi} - \frac{T}{2\pi}$ has naturally attracted his attention. He wrote me saying that he has not been able to prove the claim and asked my help. You will not find it surprising that I suggested he write to you, a suggestion he followed. Mittag-Leffler could not go beyond a point; he cannot see how you established that the series $1 - \frac{1}{2^s} - \frac{1}{3^s} - \frac{1}{5^s} + \frac{1}{6^s} - \cdots$ converges when the real part of $s$ is greater than $\frac{1}{2}$. On his behalf, I beg you to write him and free him from this embarrassment. You would him much good if you send him the clarifications, which he awaits, as he says, with impatience …

*Letter 118; Stieltjes to Hermite; Tolouse, March 30 1887.*

Sir, I had already received Mittag-Leffler's letter concerning the proposition on the number of roots of the equation $\xi(t) = 0$ and I answered with my best. I believe I was sucessful at recovering a bit of the method Riemann followed, by calculating the integral

$$\int d \log \xi(x)$$

Your letter reminds me that I had applied in my studies the formula by Weierstrass for the function $e^x - C$. For the sake of simplicity I took the constant $C$ in its form $e^a$ and, with the decomposition $\frac{1}{e^{x-a}} - 1 = \frac{e^a}{e^x - e^a}$ as auxiliary, I obtained

$$e^x - e^a = (1 - e^a) e^{-\frac{x}{e^a - 1}} \prod_{-\infty}^{+\infty} \left(1 - \frac{x}{a + 2n\pi i}\right) e^{-\frac{x}{a + 2n\pi i}}$$

For $a = 0$ there is a slight change of analytic form and we have a formula that does not differ essentially from the formula that gives $\sin x$ . Please exclude the brevity; you ought to know that I am in my final exams for the bachelor's degree.

Yours devotedly, Stieltjes.

*4.5 Stieltjes to Hermite, March 4 1887.*

*Tolouse, March 4 1887.*

Dear Sir, One thousand thanks for your letter and the translation of the Journal in Swedish. I shall return it in a few days. I especially have to ask you to consider as not valid the sending of booklets corrected by Halphen. Mr Jordan, indeed, has just sent me the manuscript and the previous folios (?). Last week, I had neither of these; in these conditions, a serious correction making sure of the exactness of all the formulae would be impossible. I have completely forgotten the meaning of certain symbols. (!) In the next few days I shall thus send you (the new these sheets; at this time I will be able to guarantee the exactness of the text

I should add a few words on Riemann's formula. ... Everyone would be morally *convinced* that $li\ (x)$ is the asymptotic expression of $F(x)$ or better of

$$G(x) = F(x) + \frac{1}{2} F\left(x^{\frac{1}{2}}\right) + \frac{1}{3} F\left(x^{\frac{1}{3}}\right) + \cdots$$

But if you ask me if this is a *proof*, in the rigorous sense of the word, the answer is no, nothing is proved in this manner.

Furthermore, to moderate the degree of trust I may have given these calculations, I would add that the agreement observed is almost too great, because one can prove rigorously that the relation

$$(A) = \frac{G(x) - li\ (x)}{x^s}$$

for $s$ a fixed number, *smaller than* $\frac{1}{2}$ (but being able to differ from $\frac{1}{2}$ by as much as one wishes), *cannot remain finite* for $x = \infty$, but must necessarily increase *beyond every limit* . Now, following the above calculations, we cannot doubt that one could believe relation remains very small or similarly tends to zero. However, it is not true that I have said is absolutely certain. The oscillating terms of Riemann's infinite series are each separate from the order $x^{\frac{1}{2}}$

If one *admits* that the *sum* of the series is similarly not of a superior order, it follows that relationship $(A)$ tends to zero when $s > \frac{1}{2}$. The calculations give much of verisimilitude to this supposition, in fact, I have said already that the accord is such that one will not doubt that the relation must increase beyond every limit, from that s is a bit smaller than $\frac{1}{2}$ . But it has not been *proven* that

$$\lim(A) = 0 \text{ when } s > \frac{1}{2}.$$

In 1885 I could not continue until having seen that certainly

$$\lim(A) = 0 \text{ when } s > \frac{3}{4}.$$

I have searched in vain to lower the limit $\frac{3}{4}$ to $\frac{1}{2}$ but I haven't doubted for a moment (on the basis of numerical calculations) that such reduction is only possible because, in fact, we have that

$$\lim(A) = 0 \text{ for } s > \frac{1}{2}.$$

Right now, without having returned to research, it seems to me that my assuredness then may have been exaggerated. In the end, the reasons offered for the limit $\frac{1}{2}$ are not serious. I freely admit that my limit ¾ is not the true limit, we know nothing about it, it is perhaps incommensurable? Here we have a great mystery. I hope one day I will learn a bit more…

---

[1] Not much is known about Stieltjes. Biographical material can be found primarily in Stieltjes 1999, but also in Rockmore 2007 and Derbyshire 2004.

[2] See *Comptes Rendus*, *Sur une function uniforme*, t.CI, July 13 1885, p. 153.

[3] See *Comptes Rendus*, Note by M. Hermite, t. CI July 13 1885, p. 112.